\documentclass[reqno,centertags,12pt]{amsart}
\usepackage{amsmath,amsthm,amscd,amssymb,latexsym,verbatim}

\usepackage{graphicx,epsf,cite}

-\marginparwidth \oddsidemargin -4mm \evensidemargin \oddsidemargin

\newtheorem{theorem}{Theorem}[section]

\theoremstyle{definition}

\theoremstyle{remark}

\newcounter{smalllist}

\allowdisplaybreaks
\numberwithin{equation}{section}

\newcommand{\lb}{\label}

\newcommand{\beq}{\begin{equation}}
\newcommand{\eeq}{\end{equation}}

\newcommand{\bal}{\begin{align}}
\newcommand{\eal}{\end{align}}
\newcommand{\bals}{\begin{align*}}
\newcommand{\eals}{\end{align*}}

\newcommand{\bbR}{{\mathbb{R}}}

\newcommand{\eps}{\varepsilon}

\begin{document}
\title[Blow up for 2D Euler Equation]
{Blow up for the 2D Euler Equation \\ on Some Bounded Domains}

\author{Alexander Kiselev}
\address{\noindent Department of Mathematics \\ Rice University \\
Houston, TX 77005, USA \newline Email: kiselev@rice.edu}

\author{Andrej Zlato\v s}
\address{\noindent Department of Mathematics \\ University of
Wisconsin \\ Madison, WI 53706, USA \newline Email:
zlatos@math.wisc.edu}


\begin{abstract}
We find a smooth solution of the 2D Euler equation on a bounded domain which exists and is unique in a natural class locally in time, but blows up in finite time in the sense of its vorticity losing continuity. The domain's boundary is smooth except at two points, which are interior cusps.
\end{abstract}

\maketitle

\section{Introduction} \lb{S1}

Consider the 2D Euler equation on a bounded domain $D\subseteq\bbR^2$, in vorticity formulation:
\beq\lb{1.1}
\omega_t+u\cdot\nabla\omega = 0, \qquad \omega(0,\cdot)=\omega_0.
\eeq
The incompressible velocity $u$ satisfies the no-flow condition $u\cdot n=0$ on $\partial D$ (with $n$ the unit outer normal to $\partial D$), and is found from the vorticity $\omega$ via
the Biot-Savart law
\beq\lb{1.2}
u(t,x)= \int_D K_D(x,y) \omega(t,y)dy,
\eeq
where the Kernel $K_D$ is given by $K_D(x,y):=\nabla^\perp G_D(x,y)$, with $\nabla^\perp=(-\partial_{x_2},  \partial_{x_1})$ and $G_D\ge 0$ the Green's function for $D$.

If the domain $D$ and the initial data $\omega_0$ are smooth, results on existence, uniqueness, and global regularity of solutions of the
2D Euler equation go back to the 1930s works of Wolibner \cite{Wol} and H\"older \cite{Hold}. Even though the solutions in this case
are globally regular, infinite-time growth of their derivatives and small scale creation appear to be ubiquitous in two dimensional fluids. 

The best known upper bound on the growth of the gradient of the vorticity is double-exponential in time (see \cite{YudDE}, even though the result is also implicit in earlier works).
The question of how fast it can actually grow has received much attention in the literature. First results of this type were due
to Yudovich \cite{Yud1,Yud2}, who constructed examples with unbounded growth of the vorticity gradient on the boundary of the domain. Recently, such growth at
the boundary was shown to be fairly generic in \cite{MSY} and \cite{Koch}. Better quantitative estimates were obtained by Nadirashvili \cite{Nad1}, who built
examples with at least  linear growth, and Denisov \cite{Den1}, whose construction gives a super-linear lower bound for the vorticity gradient. Denisov also constructed
examples with a double-exponential growth for an arbitrary but finite time \cite{Den2}. In these examples, growth happens 
in the bulk of the fluid. Very recently, Kiselev and \v Sver\' ak \cite{KS} constructed examples with a double-exponential in time growth of the vorticity gradient  
on the boundary, which is of course  the fastest possible rate of growth. Whether such growth can happen in the bulk remains open. The best result here so far is due to Zlato\v s \cite{Zlatos}, who built examples with at least exponential in time growth.

The 2D Euler equation in less regular domains and with less regular initial data has also been  studied extensively. 
Yudovich theory \cite{Yudth} (see also \cite{MB,MP}) guarantees existence and uniqueness
of a bounded weak solution $\omega$ if the domain is $C^\infty$ and the initial data $\omega_0$ is $L^\infty.$  
Uniqueness comes from the fact that the fluid velocity corresponding to a bounded vorticity $\omega$ is log-Lipshitz.
This makes bounded vorticities a very natural class  to consider (although incremental generalizations are possible). 
The smoothness assumptions on the boundary of the domain can also be weakened, and Yudovich theory can be extended to domains which are only $C^{1,1}$ \cite{GL}. 
Moreover, if the initial data and domain have a higher regularity, the solution inherits it as well \cite{MP}.
The 2D Euler equation has also been studied on certain types of domains with singular points, including domains with corners. 
We refer to recent works \cite{GL,Lacave,LMW} for some results on
existence, uniqueness, and regularity of solutions, as well as for further references.

Our main purpose in this  note is to provide an example showing that if the boundary of the domain is not {\it everywhere} sufficiently regular, 
then much more dramatic effects than just a fast growth of the derivatives of solutions can take place.
We show that there are solutions of the 2D Euler equation on domains with two cusps (and smooth elsewhere)
which are smooth locally in time but blow up in finite time, in the sense of the vorticity losing continuity. We note that our construction will
rely only on the above existence and uniqueness results for smooth domains and initial data, due to a special structure of our example.


The example involves a ``double stadium" domain. Let 
\[
\tilde{D}_r^\pm:= \left[ (-r,r)\times(\pm1-r,\pm1+r) \right] \cup B_r(-r,\pm1) \cup B_r(r,\pm1)
\]
be the ``stadium'' domain in $\bbR^2$ with width $r>0$ and centered at $(0,\pm 1)$.
While $\tilde{D}_r^\pm$ are defined explicitly, they are only $C^{1,1}.$
Let us therefore define domains $D_r^\pm$ by making the connections between the circular arcs and the horizontal intervals (which form $\partial \tilde D_r^\pm$) infinitely smooth. This can be done by adjusting the circular arcs 
near the points of contact with the intervals, and such an adjustment can be made as small as we want in $C^1.$ 
Let us do this so that  again all the $D_r^\pm$ are rescaled copies of each other, centered at $(0,\pm 1)$. 
Finally, let
\[
D':=D_1^+\cup D_1^- \cup \left[ (-1,1)\times\{0\} \right]
\]
and denote $D:=D_1^+$ the upper half of $D'$. 
%

If $\omega_0\in L^\infty(D')$ is odd in $x_2$, then the Yudovich theory shows that there is a unique odd-in-$x_2$ weak solution $\omega$ to \eqref{1.1} on $D'\times\bbR$, namely the one whose restriction to $D\times\bbR$ solves \eqref{1.1} there. Moreover, if $\omega_0 \in C(D'),$ then  it follows from Yudovich theory that $\omega(t,\cdot)\in C(D)$ 
 for all $t>0.$
The following theorem provides a necessary and sufficient condition for $\omega$ to remain continuous on the whole domain $D'.$ 

\begin{theorem}\lb{T.1.1}
Let $\omega_0\in C(\bar D')$ be non-negative on $D$ and odd in $x_2$.  Then the unique odd-in-$x_2$ solution $\omega$ to \eqref{1.1} on $D'\times\bbR$ is continuous for all $t>0$ if and only if $\omega_0(\partial D)=0$.
\end{theorem}

\it Remark. \rm If in addition we have $\omega_0 \in C^\infty (\bar D),$ then  classical regularity theory gives us $\omega(t,\cdot)\in C^\infty (D)$ for all $t>0$. 
If furthermore $\omega_0$ is supported away from $[-1,1]\times\{0\} $, then finite speed of propagation (i.e., bounded $u$) implies $\omega(t,\cdot)\in C^\infty(D')$ for all small enough $t>0$. By Theorem~\ref{T.1.1}, such a solution will become discontinuous in finite time unless $\omega_0$  vanishes on all of $\partial D.$ 


\section{Proof of Theorem \ref{T.1.1}} \lb{S2}

If $\omega_0(\partial D)=0$, then the result follows from Yudovich theory because both $D$ and $\partial D$ are invariant under the flow map generated by $u$.

Let us now assume that $\omega_0(x_0)>0$ for some $x_0\in\partial D$.
We will show that if $X'(t;x)=u(t,X(t;x))$ and $X(0;x)=x$ for $x\in\bar D$, then there is $T>0$ such that $X(T;x_0)=0$.  This, the continuity of $\omega_0$ and of $X(t;\cdot)$  on $\bar D$, and oddness of $\omega$ in $x_2$ then show that $\omega(T,\cdot)$ is discontinuous at the origin.

Let $\delta>0$ be such that  $|\omega_0^{-1}([\delta,\infty))|\ge 30\delta\, (> (4+2\pi)2\delta+\delta)$.  Because $u$ is incompressible and $\omega(t,\cdot)\ge 0$ on $D$, it follows that for each $t\ge 0$,
\beq\lb{2.1}
|\{ y\in D_{1-2\delta}^+ \,|\, \omega(t,y)\ge\delta \}| \ge \delta.
\eeq

Next let
\[
\kappa:=\inf_{x\in \partial D_{1-\delta}^+ \,\&\, y\in \bar D_{1-2\delta}^+} G_D(x,y).
\]
Continuity of $G_D$ away from the diagonal $x=y$, positivity of $G_D$ on $D$, and compactness of $\bar D_{r}^+$ show that $\kappa>0$.  Let $v>0$ solve $\Delta v=0$ on $D\setminus \bar D_{1-\delta}^+$ with $v=0$ on $\partial D$ and $v=\kappa$ on $\partial D_{1-\delta}^+$.  By the Hopf Lemma it follows that $\eps:=\inf_{x\in\partial D} |\nabla v(x)|>0$.

 We then have $G_D(x,y)\ge v(x)$ for $x\in D\setminus \bar D_{1-\delta}^+$ and $y\in  \bar D_{1-2\delta}^+$, because the inequality holds for $x\in \partial(D\setminus \bar D_{1-\delta}^+)$ (and both functions are harmonic in $x$).  So we have $|K_D(x,y)|\ge\eps$ for $x\in\partial D$ and $y\in \bar D_{1-2\delta}^+$.  This and \eqref{2.1} show that for $x\in \partial D$,  the right-hand side of \eqref{1.2} equals $a(t,x)(n_2(x),-n_1(x))$, with  $a(t,x)\ge \eps\delta^2$.  Thus $X(t;x_0)$ travels clockwise along $\partial D$ with velocity $\ge \eps\delta^2$, which means that $X(T;x_0)=0$ for some $T>0$.  The proof is finished.

\medskip

{\bf Acknowledgements.}
AK acknowledges partial support by NSF  grants DMS-1104415 and DMS-1159133, and by a Guggenheim Fellowship.
AZ acknowledges partial support  by NSF grants DMS-1056327  and DMS-1159133. 



\end{document}